# Strong subgroup chains and the Baer-Specker group

Oren Kolman[1]

**Abstract.** Examples are given of non-elementary properties that are preserved under *C*-filtrations for various classes *C* of abelian groups. The Baer-Specker group $\mathbb{Z}^\omega$ is never the union of a chain $\langle A_\alpha : \alpha < \delta \rangle$ of proper subgroups such that $\mathbb{Z}^\omega/A_\alpha$ is cotorsionfree. Cotorsionfree groups form an abstract elementary class (AEC). The Kaplansky invariants of $\mathbb{Z}^\omega/\mathbb{Z}^{(\omega)}$ are used to determine the AECs $^\perp(\mathbb{Z}^\omega/\mathbb{Z}^{(\omega)})$ and $^\perp(B/A)$, where $B/A$ is obtained by factoring the Baer-Specker group $B$ of a ZFC extension by the Baer-Specker group $A$ of the ground model, under various hypotheses, yielding information about its stability spectrum.

## Introduction

Unions of chains of subgroups of the Baer-Specker group $\mathbb{Z}^\omega$, the product of $\omega$ many copies of the infinite cyclic group $\mathbb{Z}$, have been studied in recent work of Blass and Irwin [BlIr] and Fuchs and Göbel [FG]. In [FG], the authors show that if $\langle A_\alpha : \alpha < \delta \rangle$ is a continuous ascending chain of abelian groups such that $A_0 = 0$ and for all $\alpha < \delta$, $A_{\alpha+1}/A_\alpha$ is slender, then the union $\cup_{\alpha < \delta} A_\alpha$ is slender; they deduce that the Baer-Specker group $\mathbb{Z}^\omega$ is not the union of any continuous ascending chain $\langle A_\alpha : \alpha < \delta \rangle$ of proper subgroups such that $\mathbb{Z}^\omega/A_\alpha$ is cotorsionfree, provided that $\delta < cov(\mathbf{B})$. Recall that the covering number $cov(\mathbf{B})$ is the least cardinal $\kappa$ such that $\kappa$ meagre sets cover the real line. Fuchs and Göbel also prove that $\mathbb{Z}^\omega$ is not the union of any countable ascending chain $\langle A_n : n < \omega \rangle$ of slender pure subgroups. Blass and Irwin [BlIr] prove that $\mathbb{Z}^\omega$ is never the union of a chain of proper subgroups each isomorphic to $\mathbb{Z}^\omega$, if the chain has length less than *cov*(**B**).

---

[1] I thank the referee for two detailed reports containing fundamental corrections and invaluable comments that led to a total reworking of this paper and a fuller exploitation of the pure-injectivity of $\mathbb{Z}^\omega/\mathbb{Z}^{(\omega)}$.

Eliminating the hypothesis that $\tau < cov($**B**$)$ from the Fuchs-Göbel corollary, we prove that $\mathbb{Z}^\omega$ is never the union of a chain $\langle A_\alpha : \alpha < \delta \rangle$ of proper subgroups such that $\mathbb{Z}^\omega / A_\alpha$ is cotorsionfree. The analogous question concerning the higher Baer-Specker group $\mathbb{Z}^\kappa$ for an uncountable cardinal $\kappa$ is examined briefly under additional set-theoretic axioms. Using results of Dugas and Göbel [DG], we deduce for example that if the axiom of constructibility V = L holds, then $\mathbb{Z}^\kappa$ is not the union of any continuous ascending chain $\langle A_\alpha : \alpha < \delta \rangle$ of proper subgroups such that for all $\alpha < \delta$, (1) $A_\alpha$ does not contain an isomorphic copy of $\mathbb{Z}^\kappa$ and (2) $\mathbb{Z}^\kappa / A_\alpha$ is a product. Throughout the paper, it is tacitly assumed that all groups under consideration are abelian.

Let *C* be a class of groups that is closed under isomorphism. A *C–filtration* of a group *G* is a continuous increasing chain $\langle A_\alpha : \alpha \leq \delta \rangle$ of its subgroups such that $A_0 = 0$, $A_\delta = G$, and for all $\alpha < \delta$, $A_{\alpha+1}/A_\alpha$ is isomorphic to a member of *C*. In the terminology of *C*-filtrations, the Fuchs-Göbel theorem states that if *S* is the class of slender groups and $\langle A_\alpha : \alpha \leq \delta \rangle$ is an *S*-filtration, then $A_\delta$ is slender. If a class *C* contains 0 and is closed under extensions, then any class that is also closed under *C*-filtrations satisfies some of the characteristic axioms of an *abstract elementary class* (AEC) ([Sh88] and [Gr]). In the area of module theory, AECs have been used recently as a unifying framework for classes of modules that do not necessarily possess a first-order axiomatization (see [BET] and references therein). Baldwin [Ba] asks for examples of AECs that are not given syntactically. We show that the class *K* of cotorsionfree groups is an AEC under a suitable strong submodel relation. However, *K* is never of the form $^\perp C = \{G : \text{Ext}(G, X) = 0$ for all $X \in C\}$, for any class *C*.

The paper concludes with some brief remarks on the algebraic structure of the quotient groups obtained by factoring the Baer-Specker group *B* of an extension (in the set-theoretic sense) by the Baer-Specker group *A* of the ground model. If *M* and *N* are transitive models of ZFC (Zermelo Fraenkel set theory with the axiom of choice), and *N* is an extension of *M*, then $A = (\mathbb{Z}^\omega)^M$, the Baer-Specker group in *M*, is a subgroup of $B = (\mathbb{Z}^\omega)^N$, the Baer-Specker group in *N*. The Baer-Specker quotient group, *B*/*A*, is a torsionfree cotorsion group (in *N*) and hence is pure-injective (algebraically compact). The Kaplansky invariants of *B*/*A* determine its algebraic structure completely. The class $^\perp(B/A)$ is an example of an AEC whose stability properties can be altered by forcing. We compute $^\perp(B/A)$ under various hypotheses. Applications include $^\perp(\mathbb{Z}^\omega/(\mathbb{Z}^\omega \cap L))$ when there are only countably many constructible reals, and when V = L. In the case where *N* is a forcing extension of *M*, the algebraic properties of the quotient provide a measure of the distance between the reals in the ground model and the reals in the generic extension. Profound works in this area abound: see [Bart], [RoSh] and the papers of many other authors in this tradition. We make just the simple remark that if a forcing **P** adds reals, then in *M*[*G*], there is no **K**-filtration of *B* from *A*, where **K** is the class of cotorsionfree groups in *M*[*G*].

The notation is standard and follows [EM], [K] and [BET].

Let us recall some classical results to which repeated appeals will be made throughout the paper. We shall call a group *G* a *product* if *G* is isomorphic to $\mathbb{Z}^\lambda$ for some cardinal $\lambda$. A group *G* is a *countable product* if *G* is isomorphic to $\mathbb{Z}^X$ for some $X \in \omega \cup \{\omega\}$. The following important results, due to Nunke, Dugas and Göbel, and Łoś, will prove invaluable; they are mostly collected in [EM], [F1] and [F2].

**Theorem** (Nunke [N2])
Every epimorphic image of the Baer-Specker group is a direct sum of a cotorsion group and a countable product.

**Theorem** (Nunke [N1]; Dugas, Göbel [DG])
Suppose that $\kappa$ is not $\omega$-measurable and *A* is a subgroup of $\mathbb{Z}^\kappa$ such that $\mathbb{Z}^\kappa/A$ is a product. Then *A* is a direct summand of $\mathbb{Z}^\kappa$ if any of the following conditions holds:
1. $\kappa = \omega$;
2. every Whitehead group of cardinality at most $\kappa$ is free;
3. $\mathbb{Z}^\kappa/A$ is a countable product.

Note that condition (2) subsumes condition (1) since countable Whitehead groups are free.

**Theorem** (Łoś; see [F2])
Suppose that $\kappa$ is an infinite cardinal that is not $\omega$-measurable. No proper direct summand of $\mathbb{Z}^\kappa$ contains $\mathbb{Z}^{(\kappa)}$.

# 1 Filtrations by cotorsionfree groups

**Definition**
Let $C$ be a class of groups that is closed under isomorphism, and suppose that $A$ is a subgroup of $B$.
(1) A *C–filtration* of $B$ from $A$ is a continuous increasing chain $\langle A_\alpha : \alpha \leq \delta \rangle$ of subgroups of $B$ such that $A_0 = A$, $A_\delta = B$, and for all $\alpha < \delta$, $A_{\alpha+1}/A_\alpha$ is isomorphic to a member of $C$.
(2) $B$ is *C–filtrable* from $A$ if there is a *C–filtration* of $B$ from $A$.

Thus, a *C–filtration* of $B$ is a *C–filtration* of $B$ from 0. Also if $\langle A_\alpha : \alpha \leq \delta \rangle$ is a *C–filtration* of $B$ from $A$, then $\langle A_\alpha/A : \alpha \leq \delta \rangle$ is a *C–filtration* of $B/A$.

**Definition**
A class $C$ of groups is closed under extensions if $G/A \in C$ and $A \in C$ imply $G \in C$.

**Proposition 1.1.**
Suppose $C$ is closed under extensions. The following are equivalent:
(1) For all $A$ and $B$, for every *C–filtration* $\langle A_\alpha : \alpha \leq \delta \rangle$ of $B$ from $A$, $A \in C$ implies $B \in C$.
(2) For all $A$ and $B$, for every *C–filtration* $\langle A_\alpha : \alpha \leq \delta \rangle$ of $B$ from $A$, $A \in C$ and $\alpha < \delta$ imply $B/A_\alpha \in C$.

Proof. For (1) $\Rightarrow$ (2), apply (1) to the *C–filtration* $\langle A_\beta/A_\alpha : \alpha \leq \beta < \delta \rangle$ of $B/A_\alpha$; for (2) $\Rightarrow$ (1), use (2) to deduce that $B/A = B/A_0 \in C$, and then recall that $C$ is closed under extensions. ♦

Recall that the class of cotorsionfree groups is closed under subgroups, extensions and direct products. In [FG], the authors observe and implicitly prove that if $K$ is the class of cotorsionfree groups and $\langle A_\alpha : \alpha \leq \delta \rangle$ is a $K$-filtration of $B$ from $A \in K$, then $B \in K$. We therefore just sketch the argument, making explicit the role that closure of the class $K$ under extensions plays, since this will be of interest in the discussion of abstract elementary classes.

**Theorem 1.2.**
Let $K$ be the class of cotorsionfree groups. If there exists a $K$-filtration $\langle A_\alpha : \alpha \leq \delta \rangle$ of $B$ from $A \in K$, then $B \in K$.

Proof. We prove simultaneously by induction on $\delta$ the following two claims:
   (1) $B \in K$, and (2) for all $\alpha < \delta$, $B/A_\alpha \in K$.
If $\delta = 0$, then the result simply restates the hypothesis that $A$ is cotorsionfree.
If $\delta = \gamma + 1$, then $B$ is cotorsionfree, being an extension of $A_\gamma$ by a cotorsionfree group and $B/A_\alpha$ is cotorsionfree since the quotient $(B/A_\alpha)/(A_\gamma/A_\alpha)$ is cotorsionfree and $K$ is closed under extensions.
If $\delta$ is a limit ordinal, then for all $\alpha < \delta$, $A_\alpha$ is a cotorsionfree pure subgroup of $B$. It is easy to see that neither $B$ nor $B/A_\alpha$ can contain $\mathbb{Q}$ or $\mathbb{Z}_p$ for any prime $p$; if, for a contradiction, $B$ contains $J_p$ for some prime $p$, note that since $A_{\alpha+1}/A_\alpha$ is cotorsionfree, then as in [FG] there exists $\alpha < \beta < \delta$ such that $A_\beta/A_\alpha$ is not cotorsionfree, contradicting the induction hypothesis for (2) applied to $A_\beta/A_\alpha$. Therefore $B$ cannot contain $J_p$; by the Proposition 1.2, $B/A_\alpha$ is also cotorsionfree. ♦

It will be convenient to have a separate statement of (2) for later.

**Corollary 1.3.**
Let $K$ be the class of cotorsionfree groups. Suppose that $\langle A_\alpha : \alpha \leq \delta \rangle$ is a $K$-filtration of $B$ from $A \in K$. Then for all $\alpha < \delta$, $B/A_\alpha \in K$. ♦

**Theorem 1.4.**
The Baer-Specker group $\mathbb{Z}^\omega$ is never the union of a chain $\langle A_\alpha : \alpha < \delta \rangle$ of proper subgroups such that $\mathbb{Z}^\omega/A_\alpha$ is cotorsionfree.

Proof. Suppose that $\langle A_\alpha : \alpha < \delta \rangle$ is a chain of proper subgroups such that $\mathbb{Z}^\omega/A_\alpha$ is cotorsionfree. If $cf(\delta) = \omega$, let $\langle \alpha_n : n < \omega \rangle$ be a cofinal sequence in $\delta$. Then by [FG], since $\omega < cov(\mathbf{B})$, $\mathbb{Z}^\omega \neq \cup_{n<\omega} A_{\alpha_n} = \cup_{\alpha<\delta} A_\alpha$. If $cf(\delta) > \omega$, then for each $n$, there exists $\alpha_n$ such that $e_n \in A_{\alpha_n}$ (recall $e_n$ is the unit vector with 1 at the $n$-th place, 0 elsewhere); so $\alpha^* = \sup \langle \alpha_n : n < \omega \rangle < \delta$, and $\mathbb{Z}^{(\omega)}$ is contained in $A_{\alpha^*}$. But $\mathbb{Z}^\omega/\mathbb{Z}^{(\omega)}$ is cotorsion and therefore so is its quotient $\mathbb{Z}^\omega/A_{\alpha^*}$, contradicting the hypothesis that $\mathbb{Z}^\omega/A_{\alpha^*}$ is cotorsionfree. ♦

In the above result, the chain is not assumed continuous.

**Corollary 1.5.**
Suppose that $\langle A_\alpha : \alpha < \delta \rangle$ is a chain of proper direct summands of $\mathbb{Z}^\omega$. Then $\mathbb{Z}^\omega \neq \cup_{\alpha<\delta} A_\alpha$. ♦

**Corollary 1.6.**
The Baer-Specker group $\mathbb{Z}^\omega$ is never the union of a chain $\langle A_\alpha : \alpha < \delta \rangle$ of slender direct summands. ♦

Note however that $\mathbb{Z}^\omega$ is the union of a continuous ascending chain of length $cf(2^{\aleph_0})$ of slender subgroups; also, by Kulikov's theorem, $\mathbb{Z}^\omega$ is the union of a ascending chain of countable length of subgroups that are free (hence slender) (see [F1]).

**Corollary 1.7.**
Let $U$ be a non-principal ultrafilter on $\omega$. Then the ultrapower $\mathbb{Z}^\omega/U$ is never the union of a chain $\langle B_\alpha : \alpha < \delta \rangle$ of proper subgroups such that $(\mathbb{Z}^\omega/U)/B_\alpha$ is cotorsionfree.

Proof. If $B_\alpha = A_\alpha/K_U$, where $A_\alpha$ is a proper subgroup of $\mathbb{Z}^\omega$ and $K_U = \{x \in \Pi_{i \in I} M_i : \{i \in I : x(i) = 0\} \in U\}$, then $\mathbb{Z}^\omega/A_\alpha \cong (\mathbb{Z}^\omega/U)/B_\alpha$ is cotorsionfree; apply Theorem 1.4 to the chain $\langle A_\alpha : \alpha < \delta \rangle$. ♦

We consider next whether these corollaries about $\mathbb{Z}^\omega$ have provable analogues in higher cardinalities. This will require information on the quotients of the higher Baer-Specker group $\mathbb{Z}^\kappa$ for an uncountable cardinal $\kappa$. Recall the well-known fact that the axiom of constructibility V = L implies that measurable (and hence $\omega$-measurable) cardinals do not exist.

**Theorem** (Dugas, Göbel [DG])
If V = L, then for all cardinals $\kappa$ and all subgroups $A$ of the higher Baer-Specker group $\mathbb{Z}^\kappa$, $A$ is a direct summand of $\mathbb{Z}^\kappa$ if and only if $\mathbb{Z}^\kappa/A$ is a product.

**Theorem 1.7.**
Assume that V = L holds. Suppose that $\langle A_\alpha : \alpha < \delta \rangle$ is a continuous ascending chain of subgroups of $\mathbb{Z}^\kappa$ such that for all $\alpha < \tau$ (1) $A_\alpha$ is not isomorphic to $\mathbb{Z}^\kappa$ and (2) $\mathbb{Z}^\kappa/A_\alpha$ is a product. Then $\mathbb{Z}^\kappa \neq \cup_{\alpha < \delta} A_\alpha$.

Proof. Case 1: $\kappa < cf(\delta)$. Then for some $\alpha^* < \delta$, $\mathbb{Z}^{(\kappa)}$ is contained in $A_{\alpha^*}$. By the results of Dugas and Göbel and Łoś, it follows that $A_{\alpha^*} = \mathbb{Z}^\kappa$, a contradiction.
Case 2: $cf(\delta) \leq \kappa$. Passing to a continuous cofinal subchain if necessary, we may assume that $cf(\delta) = \delta$. Since $\delta \leq \kappa < cf(|\mathbb{Z}^\kappa|)$, for some $\alpha^* < \delta$, $A_{\alpha^*}$ has cardinality $2^\kappa$ and is a direct summand of $\mathbb{Z}^\kappa$; referring to Theorem 1.4 on pages 294-295 in [EM] and recalling that GCH is a consequence of V = L, it follows that $A_{\alpha^*}$ must be isomorphic to $\mathbb{Z}^\kappa$, contradicting (1). ♦

The referee generously provided a direct proof of the last assertion and for convenience it is reproduced here with thanks.

Claim: if $\kappa$ is not $\omega$-measurable, then any direct summand of $\mathbb{Z}^\kappa$ must be isomorphic to $\mathbb{Z}^\lambda$ for some $\lambda \leq \kappa$. Why? Suppose $\mathbb{Z}^\kappa$ is isomorphic to $A \oplus B$. The non-$\omega$-measurability hypothesis implies that $\mathbb{Z}^\kappa$ is reflexive, and therefore so are the summands $A$ and $B$. Non-$\omega$-measurability also implies that the dual group $(\mathbb{Z}^\kappa)^*$ is free of rank $\kappa$, and therefore its summands $A^*$ and $B^*$ are also free. Let $\lambda$ be the rank of $A^*$. Then $\lambda \leq \kappa$ and

$$A \cong A^{**} \cong (\mathbb{Z}^{(\lambda)})^* \cong \mathbb{Z}^\lambda,$$

as claimed.

Since GCH follows from $V = L$, if $\lambda < \kappa$, then $|A_{\alpha^*}| = |\mathbb{Z}^\lambda| = \lambda^+ \leq \kappa < 2^\kappa$. In case (2) of the proof above, $|A_{\alpha^*}| = 2^\kappa$, so we can conclude that $\lambda = \kappa$.

**Theorem** (Dugas, Göbel [DG])
Suppose that $\kappa$ and $\lambda$ are cardinals which are not $\omega$-measurable and $A$ is a subgroup of the higher Baer-Specker group $\mathbb{Z}^\kappa$ which is a homomorphic image of $\mathbb{Z}^\lambda$. Then $\mathbb{Z}^\kappa/A$ is the direct sum of a cotorsion group and a product.

**Corollary 1.8.**
Assume that $V = L$ holds. Suppose that for all $\alpha < \delta$, $\lambda_\alpha < \kappa$, and $\langle A_\alpha : \alpha < \delta \rangle$ is a continuous ascending chain of subgroups of $\mathbb{Z}^\kappa$ such that (1) $A_\alpha$ is a homomorphic image of $\mathbb{Z}^{\lambda_\alpha}$ and (2) $\mathbb{Z}^\kappa/A_\alpha$ is cotorsionfree. Then $\mathbb{Z}^\kappa \neq \cup_{\alpha < \delta} A_\alpha$.

Proof. Since GCH holds in L and $\lambda_\alpha < \kappa$, $A_\alpha$ cannot contain an isomorphic copy of $\mathbb{Z}^\kappa$ and by (2) $\mathbb{Z}^\kappa/A_\alpha$ is a product; an appeal to Theorem 1.7 completes the proof. ♦

## 2  Abstract elementary classes and cotorsionfreeness

Let us return to the class **K** of cotorsionfree groups and the natural question whether there exists a partial order $<_K$ on **K** such that (**K**, $<_K$) is an abstract elementary class. First, recall the concept of an elementary class. Let $C$ be a class of structures all of the same similarity type L($C$). The class $C$ is an *elementary class* if there exists a first-order theory T in L($C$) such that $C$ = Mod(T), the class of models of T.

For definiteness, the similarity type of groups is $\{+, -, 0\}$; for $\mathbb{Z}$-modules (or generally $R$-modules), add a unary function symbol for each element $r$ of the ring to express scalar multiplication by $r$. Many natural classes of abelian groups are non-elementary. For example, since every countable elementary submodel of the group $J_p$ of the $p$-adic integers is cotorsionfree, it follows that the class $K$ of cotorsionfree groups is non-elementary. Similarly, by Łoś Lemma, $\mathbb{Z}$ is elementarily equivalent to the ultrapower $\mathbb{Z}^\omega/U$ by a non-principal ultrafilter $U$, and so the class of slender groups is non-elementary (in fact it is not pseudo-elementary).

We shall write $M \subseteq N$ to mean that $M$ is a substructure of $N$.

**Definition**
Let $C$ be a class of structures all of the same similarity type L($C$). The ordered pair ($C$, $<_C$) is an *abstract elementary class* (AEC) if the following axioms are satisfied:
(A0) Both $C$ and the binary relation $<_C$ on $C$ are closed under isomorphism, i.e. (i) if $f$ is an isomorphism of $M \in C$ onto $N$, then $N \in C$, and (ii) if $f_i$ ($i = 0, 1$) is an isomorphism from $M_i$ onto $N_i$, $M_0 <_C M_1$, and $f_1$ extends $f_0$ (i.e., $f_0 \subseteq f_1$), then $N_0 <_C N_1$.
(A1) $<_C$ is a partial order on $C$.
(A2) if $M <_C N$, then $M \subseteq N$.
(A3) if $\langle A_\alpha : \alpha < \delta \rangle$ is a continuous $<_C$-increasing chain, then:
    (1) $\cup_{\alpha < \delta} A_\alpha \in C$;
    (2) for each $\beta < \delta$, $A_\beta <_C \cup_{\alpha < \delta} A_\alpha$;
    (3) if for each $\beta < \delta$, $A_\beta <_C M \in C$, then $\cup_{\alpha < \delta} A_\alpha <_C M$.
(A4) If $A, B, C \in C$, $A <_C C$, $B <_C C$, and $A \subseteq B$, then $A <_C B$.
(A5) There is a Löwenheim–Skolem number LS($C$) such that if $A \subseteq B \in C$, there exists $A' \in C$ such that $A \subseteq A'$, $A' <_C B$, and $|A'| \leq |A| + $ LS($C$).

We say that $\langle A_\alpha : \alpha < \delta \rangle$ is a continuous $<_C$-increasing chain if for all $\alpha < \delta$, $A_\alpha \in C$, $A_\alpha <_C A_{\alpha+1}$, and if $\xi$ is a limit ordinal then $A_\xi = \cup_{\zeta < \xi} A_\zeta$. If $M <_C N$, then $M$ is called a *strong submodel* of $N$.

Numerous examples of AECs are given in [Gr], [BaLe] and [BET]. See [Gr] for an explanation of AECs in the context of classification theory for non-elementary classes. We recall from [BET] the following definition and notation. For a class $C$ of groups, let $(^\perp C, <_C)$ be defined as follows: $^\perp C = \{G : \mathrm{Ext}(G, X) = 0 \text{ for all } X \in C\}$, and the partial order $<_C$ is defined by $G <_C H$ if $G$ is a subgroup of $H$ and $G$ and $H/G$ belong to $^\perp C$. Baldwin, Eklof and Trlifaj ([BET], Theorem 1.20) prove that $(^\perp C, <_C)$ is an AEC if and only if every member of $C$ is cotorsion.

We seek a strong subgroup relation $<_K$ on the class $K$ of cotorsionfree groups under which $(K, <_K)$ is an abstract elementary class. Theorem 1.1 and its corollary express the fact that $K$ satisfies the axioms (A3)(1) and (A3)(2) of an abstract elementary class under the partial order defined as follows: $G < H$ if $G$ is a subgroup of $H$ and $G$ and $H/G$ belong to $K$. The axioms (A0), (A1), (A2) and (A4) are also evident for $(K, <_K)$. The first problematic axiom is (A3)(3): if we take $M = \mathbb{Z}^\omega \in K$, and identify $A_n = \mathbb{Z}^n$ with its inclusion in $\mathbb{Z}^\omega$, then $A_n < A_{n+1}$, $\cup_{n<\omega} A_n = \mathbb{Z}^{(\omega)}$, but $\mathbb{Z}^\omega / \mathbb{Z}^{(\omega)}$ is a cotorsion group.

**Theorem** ([KolSh])

There exists a sentence $\psi$ in the infinitary language $L_{\lambda, \omega}$ with $\lambda = (2^{\aleph_0})^+$ such that $K = \mathrm{Mod}(\psi)$.

Fix a fragment $\mathbb{A}$ of $L_{\lambda, \omega}$ containing the sentence $\psi$ and define $G <_K H$ if and only if $G$ is an $L_\mathbb{A}$-elementary substructure of $H$.

**Theorem 2.1.**
The class $(K, <_K)$ of cotorsionfree groups is an abstract elementary class, where $G <_K H$ if and only if $G$ is an $L_\mathbb{A}$-elementary substructure of $H$. The Löwenheim–Skolem number $LS(K)$ is $|\mathbb{A}|$.

Proof. Immediate from the quoted theorem and well-known model theory of the infinitary language $L_{\lambda, \omega}$. ♦

Since every Whitehead group is slender and hence cotorsionfree, but the cotorsionfree group $\mathbb{Z}^\omega$ is not Whitehead, it follows that $^\perp\mathbb{Z}$ is a proper subclass of $K$. However, in answer to a question of the referee, $K$ cannot be represented as $^\perp C$ for any class $C$ of abelian groups.

**Corollary 2.2.**
The class $K$ of cotorsionfree groups is never of the form $^\perp C$ for any class $C$ of abelian groups.

Proof. The class $K$ of cotorsionfree groups is closed under subgroups and products, and it contains $\mathbb{Z}$, but not $\mathbb{Z}^\omega/\mathbb{Z}^{<\omega}$. If $K = {^\perp C}$, then by Lemma 4.3.17 in [GT], $\mathbb{Z}^\omega/\mathbb{Z}^{<\omega} \in {^\perp C}$. ♦

## 3  Baer-Specker quotient groups

We wish to exhibit some examples of AECs that arise from the quotient groups obtained by factoring the Baer-Specker group $B$ of an extension (in the set-theoretic sense) by the Baer-Specker group $A$ of the ground model. The Baer-Specker quotient group, $B/A$, is a torsionfree cotorsion group (in $N$) and hence is pure-injective. By [BET], the class $^\perp(B/A)$ is an AEC in $N$. We use the Kaplansky invariants of $B/A$ to compute $^\perp(B/A)$ under various hypotheses, yielding information about the stability properties of the class. The following facts will be the cornerstone of the calculations:

**Theorem** (Hulanicki [GoHu]; Balcerzyk [Balc])
(1) The group $\mathbb{Z}^\omega/\mathbb{Z}^{(\omega)}$ is pure-injective (algebraically compact).
(2) The invariants of the group $\mathbb{Z}^\omega/\mathbb{Z}^{(\omega)}$ are $\alpha_{p,n} = 0$, $\beta_p = 2^{\aleph_0}$, $\gamma_p = 0$, $\delta = 2^{\aleph_0}$, and $\mathbb{Z}^\omega/\mathbb{Z}^{(\omega)}$ is isomorphic to $\Pi_{p \in P} A_p \oplus \mathbb{Q}^{(\delta)}$ where for each prime $p$, $A_p$ is the $p$-adic completion of $J_p^{(\beta_p)}$.

**Theorem 3.1.**
(1) $^\perp\mathbb{Q}$ is the class of all abelian groups.
(2) $^\perp(\mathbb{Z}^\omega/\mathbb{Z}^{(\omega)}) = \bigcap_{p \in P} {^\perp J_p}$.

(3) $^\perp(\mathbb{Z}^\kappa/B^{<\infty}) = {}^\perp\mathbb{Q}$, where $B^{<\infty}$ is the subgroup of bounded functions.

Proof. (1) The group $\mathbb{Q}$ is divisible, hence $\mathrm{Ext}(G, \mathbb{Q}) = 0$ for every abelian group $G$.
(2) The $p$-adic completion of $J_p^{(\beta_p)}$ is $J_p^\omega$ since $\beta_p = 2^{\aleph_0}$, and $\mathbb{Q}^{(\delta)}$ is isomorphic to $\mathbb{Q}^\omega$ since $\delta = 2^{\aleph_0}$. The result follows on applying the isomorphism $\mathrm{Ext}(G, \Pi_{i \in I} C_i) \cong \Pi_{i \in I} \mathrm{Ext}(G, C_i)$.

(3) It is a remark of P. Hill (see [Ge]) that the group $\mathbb{Z}^\kappa/B^{<\infty}$ is divisible; hence it is of form $\mathbb{Q}^{(\lambda)}$ which is divisible. ♦

Recall from [BET] that for a set $P$ of maximal ideals of $\mathbb{Z}$, $\mathbf{K}(P)$ is defined as the class of all abelian groups that are $d$-torsionfree for all $d \in P$.

**Corollary 3.2.**

(1) If CH holds, then $^\perp(\mathbb{Z}^\omega/\mathbb{Z}^{(\omega)})$ is stable in $\aleph_1$ (in fact in all $\aleph_n$, $0 < n < \omega$), but not stable in $\aleph_\omega$.

(2) If CH fails, then $^\perp(\mathbb{Z}^\omega/\mathbb{Z}^{(\omega)})$ is not stable in $\aleph_1$.

Proof. By [BET], the AEC $^\perp(\mathbb{Z}^\omega/\mathbb{Z}^{(\omega)})$ is $(\mathbf{K}(P), <_{K(P)})$ for some set $P$ of maximal ideals, and the latter is stable in an infinite cardinal $\lambda$ if and only if $\lambda^{\aleph_0} = \lambda$. ♦

Suppose that $M = (M, \in^M)$ and $N = (N, \in^N)$ are models of ZF and that $N$ is an extension of $M$, i.e., every element of $M$ belongs to $N$ and $\in^N$ of $N$ agrees on $M \times M$ with $\in^M$. There are thus two Baer-Specker groups, $A = (\mathbb{Z}^\omega)^M$, the Baer-Specker group in $M$, and $B = (\mathbb{Z}^\omega)^N$, the Baer-Specker group in $N$. We shall also suppose that $M$ and $N$ are transitive models of ZF. Note that if $M$ is a transitive model of ZF, then $(\mathbb{Z}^\omega)^M = \mathbb{Z}^\omega \cap M$, where $\mathbb{Z}^\omega$ is the real Baer-Specker group, i.e., $(\mathbb{Z}^\omega)^V$. By absoluteness, $(\mathbb{Z}^{(\omega)})^M = (\mathbb{Z}^{(\omega)})^N = \mathbb{Z}^{(\omega)}$, and hence $\mathbb{Z}^{(\omega)}$ is always a subgroup of $A$. Let us call $B/A$ the *Baer-Specker quotient group*.

To determine the structure of $B/A$, we shall appeal to some well-known classical results.

**Theorem** (Chase [Ch])

Suppose that $H$ is a countable pure subgroup of $\mathbb{Z}^\omega$ and $H$ is dense in the product topology, where $\mathbb{Z}$ is equipped with the discrete topology. Then there exists an automorphism $\alpha \in Aut(\mathbb{Z}^\omega)$ such that $\alpha$ maps $H$ onto $\mathbb{Z}^{(\omega)}$.

**Corollary 3.3.**

Suppose that $H$ is a countable pure subgroup of $\mathbb{Z}^\omega$ and $H$ is dense in the product topology. Then:

(1) $\mathbb{Z}^\omega/H$ is isomorphic to $\Pi_{p \in P} A_p \oplus \mathbb{Q}^{(\delta)}$ where for each prime $p$, $A_p$ is the completion of $J_p^{(\beta_p)}$, $\beta_p = 2^{\aleph_0}$, and $\delta = 2^{\aleph_0}$.

(2) $^\perp(\mathbb{Z}^\omega/H) = \cap_{p \in P}{}^\perp J_p$.

Proof. By the theorem of Chase, $\mathbb{Z}^\omega/H$ is isomorphic to the cotorsion group $\mathbb{Z}^\omega/\mathbb{Z}^{(\omega)}$; now apply the result of Balcerzyk and Corollary 3.1. ♦

Observe that $A$ is a subgroup of $B$, and hence in $N$, $B/A$ is a cotorsion group. Since $A = \mathbb{Z}^\omega \cap M$, it follows that $A$ is pure in $B$ and hence $B/A$ is torsionfree, so that $B/A$ is algebraically compact (pure-injective). Kaplansky's classification theorem for algebraically compact groups (see [EM] or [F1]) applied in $N$ says that in $N$, $B/A$ has the form $\Pi_{p \in P} A_p \oplus \mathbb{Q}^{(\delta)}$, where for each prime $p$, $A_p$ is the completion of $J_p^{(\beta_p)}$, $0 \le \beta_p \le 2^{\aleph_0}$ for each prime $p$, and $0 \le \delta \le 2^{\aleph_0}$.

**Proposition 3.4.**
Suppose that $M$ is a transitive model of ZF. Then the Baer-Specker quotient group $B/A$ is pure-injective in $N$. If $A$ is countable, then the cardinal invariants of $B/A$ are $\alpha_{p,n} = 0$, $\beta_p = 2^{\aleph_0}$, $\gamma_p = 0$, $\delta = 2^{\aleph_0}$. $^\perp(B/A) = \cap_{p \in P}{}^\perp J_p$.

Proof. $A/\mathbb{Z}^{(\omega)}$ is a pure subgroup of the cotorsion group $B/\mathbb{Z}^{(\omega)}$ in $N$. ♦

Now consider what happens when $N$ is obtained as a generic extension of $M$ by some forcing **P**.

**Corollary 3.5.**
Suppose that $M$ is a transitive model of ZF + CH. Suppose that the forcing **P**

collapses $\aleph_1$. Then in $M[G]$, $^\perp(B/A) = \bigcap_{p \in P} {}^\perp J_p$. The stability spectrum of $^\perp(B/A)$ can be calculated from the cardinal arithmetic of $M[G]$. ♦

If the forcing **P** adds reals, then the Baer-Specker group acquires new elements in the generic extension $M[G]$, and the Baer-Specker group $A$ of the ground model may become slender. If the forcing is iterated, then new reals may appear at the stages of the iteration, giving rise to an ascending chain of subgroups of $B$ in the generic extension. A final observation about $B/A$ shows that if **P** adds reals, then in $M[G]$, $B$ is never **K**-filtrable from $A$, where **K** is the class of cotorsionfree groups in $M[G]$.

**Proposition 3.6.**
Supose that the forcing **P** adds reals. Then in $M[G]$, there is no **K**-filtration of $B$ from $A$, where **K** is the class of cotorsionfree groups in $M[G]$.

Proof. Suppose towards a contradiction that in $M[G]$ $\langle A_\alpha : \alpha \leq \delta \rangle$ is a **K**-filtration of $B$ from $A$. Then $A \in \mathbf{K}$, and hence by Corollary 1.3, $B/A \in \mathbf{K}$. But $B/A$ is cotorsion, hence 0. ♦

Let us next suppose that $M = L$ and $N = V$, i.e. $M$ is the universe of constructible sets and $N$ is the real universe. We shall call the group $\mathbb{Z}^\omega \cap L$ the constructible Baer-Specker group. Recall that reduced torsionfree groups of cardinality less than continuum are slender (Sąsiada [Są]), as are direct sums of slender groups (Fuchs [F2]); also cotorsionfree groups are reduced and torsionfree; see [EM]. The following proposition records some basic facts about the constructible Baer-Specker group. For more detailed information on L and $0^\#$, see [D], [M], [K], [EM] or [J]; $\aleph_1^{(L)}$ denotes the first uncountable cardinal in L.

**Proposition 3.7.**
1. The constructible Baer-Specker group is a pure subgroup of $\mathbb{Z}^\omega$ properly containing $\mathbb{Z}^{(\omega)}$.
2. If $V = L$, then the constructible Baer-Specker group is not slender, and hence not free.
3. If $\aleph_1^L = \aleph_1$ and $2^{\aleph_0} > \aleph_1$, then the constructible Baer-Specker group is

slender, but not free. If $\aleph_1^L < \aleph_1$, then the constructible Baer-Specker group is free.

4. If there exists a non-constructible real, then $\mathbb{Z}^\omega \cap L$ is not a direct summand of $\mathbb{Z}^\omega$.
5. If $0^\#$ exists, then the constructible Baer-Specker group is free.
6. If there exists a measurable cardinal, then the constructible Baer-Specker group is free.

Proof. Most of these are either trivial or evident from well-known results concerning L and $0^\#$.

For (3), recall that $\mathbb{Z}^\omega \cap L$ is cotorsionfree and, having cardinality $\aleph_1 < 2^{\aleph_0}$, is therefore slender. The subgroup $H$ of $\mathbb{Z}^\omega \cap L$ consisting of elements whose tails are divisible by arbitrarily large powers of a fixed prime $p$ is uncountable and non-free. The second part of (3) is immediate on recalling that for $\alpha = \aleph_1^L$, $\mathbb{Z}^\omega \cap L \subseteq L_{\alpha'}$ so if $\alpha = \aleph_1^L < \aleph_1$, then $\mathbb{Z}^\omega \cap L$ is countable and free, since $\mathbb{Z}^\omega$ is $\aleph_1$-free.

The claim (4) is immediate from (1) and the theorem of Łoś quoted in the introduction.

Since $\mathbb{Z}^\omega$ is $\aleph_1$-free and the hypotheses of (5) and (6) separately imply that $\mathbb{Z}^\omega \cap L$ is countable (in V), it follows that $\mathbb{Z}^\omega \cap L$ is free. ♦

The invariants of the quotient $\mathbb{Z}^\omega/(\mathbb{Z}^\omega \cap L)$ are easily seen to be independent of ordinary set theory. If $V = L$, then all the invariants are zero and by a result of Baldwin et al. [BCGVW-A] quoted in [BET], $\mathbf{K}(\varnothing) = {}^\perp(\mathbb{Z}^\omega/(\mathbb{Z}^\omega \cap L))$ is the class of all abelian groups and is stable. If $\aleph_1^L < \aleph_1$, then the cardinal invariants of $\mathbb{Z}^\omega/(\mathbb{Z}^\omega \cap L)$ are $\alpha_{p,n} = 0$, $\beta_p = 2^{\aleph_0}$, $\gamma_p = 0$, $\delta = 2^{\aleph_0}$; ${}^\perp(\mathbb{Z}^\omega/(\mathbb{Z}^\omega \cap L))$ is $\cap_{p \in P} {}^\perp J_p$, and the stability spectrum depends on the cardinal arithmetic of V.

**Author information**

Oren Kolman
Laboratoire de Mathématiques Nicolas Oresme - CNRS
Université de Caen BP 5186
14032 Caen Cedex, France.
E-mail: okolman@member.ams.org